\documentclass[12pt]{amsart}

\usepackage{amssymb}
\usepackage{mathrsfs}
\newcounter{fact_count}

\newcounter{conj_count}

\newtheorem{thm}{Theorem}[section]
\newtheorem{lem}[thm]{Lemma}
\newtheorem{example}[thm]{Example}
\newtheorem{prop}[thm]{Proposition}
\newtheorem{fact}[fact_count]{Fact}
\newtheorem{claim}[thm]{Claim}
\newtheorem{conjecture}[conj_count]{Conjecture}
\newtheorem{conjecture-section}[thm]{Conjecture}

\newtheorem{cor}[thm]{Corollary}
\newtheorem*{thm*}{Theorem}
\newtheorem*{defn}{Definition}

\newtheorem{question}[thm]{Question}

\theoremstyle{remark}
\newtheorem*{remark}{Remark}

\renewcommand\mathcal{\mathscr}
\newcommand\Acal{\mathcal{A}}
\newcommand\Bcal{\mathcal{B}}
\newcommand\Ccal{\mathcal{C}}

\newcommand\Fcal{\mathcal{F}}

\newcommand\Rcal{\mathcal{R}}

\newcommand\Ucal{\mathcal{U}}
\newcommand\Xcal{\mathcal{X}}

\newcommand\Rbb{\mathbb{R}}

\newcommand{\dom}{\operatorname{dom}}
\newcommand{\ran}{\operatorname{ran}}

\newcommand{\osc}{\operatorname{osc}}
\newcommand{\Osc}{\operatorname{Osc}}

\newcommand{\tr}{\operatorname{Tr}}
\newcommand{\meet}{\wedge}
\newcommand{\bigmeet}{\bigwedge}
\newcommand{\join}{\vee}

\newcommand{\Th}{{}^{\textrm{th}}}
\newcommand\axiom{\mathrm}
\newcommand\MA{\axiom{MA}}
\newcommand\PFA{\axiom{PFA}}

\newcommand\ZFC{\axiom{ZFC}}

\renewcommand{\>}{\rangle}
\renewcommand\rho{\varrho}
\newcommand\mand{\land}

\newcommand\fin{[\omega_1]^{<\aleph_0}}
\newcommand\ww{{\omega}^\omega}
\newcommand\seq{\omega^{<\omega}}
\renewcommand\epsilon{\varepsilon}
\renewcommand\mod{\operatorname{mod}}
\renewcommand\diamond{\diamondsuit}

\newcommand\typezero{\big} 
\newcommand\typeone{\big} 
\newcommand\typetwo{\big} 
\newcommand\typethree{\big} 
\newcommand\typefour{\big} 

\title[L space]{A solution to the L space problem and related ZFC 
constructions}

\author{Justin Tatch Moore}

\address{Department of Mathematics \\
Boise State University \\
Boise, ID 83724}

\email{justin@math.boisestate.edu}


\thanks{This paper is dedicated to Stevo Todor\v{c}evi\'{c} for teaching
me how to traverse $\omega_1$ and for his inspirational \cite{basis_problems}.
I would like to thank Boban Veli\v{c}kovi\'c for his careful reading
of early drafts of this paper.
The clarity of the present arguments owe much to Todor\v{c}evi\'{c}'s and
Veli\v{c}kovi\'{c}'s suggestions.
A number of other people --- too many to name --- have also made the
generous contribution of reading the paper and offering suggestions.
The research presented in this paper was funded by NSF grant
DMS--0401893.  
}
\begin{document}

\begin{abstract}
In this paper I will construct a non-separable
hereditarily Lindel\"of space (L space) without any additional
axiomatic assumptions.
I will also show that there is a
function $f:[\omega_1]^2 \to \omega_1$ such that
if $A,B \subseteq \omega_1$ are uncountable and $\xi < \omega_1$,
then there are $\alpha < \beta$ in $A$ and $B$ respectively with
$f (\alpha,\beta) = \xi$.
Previously it was unknown whether such a function existed even if
$\omega_1$ was replaced by $2$.
Finally, I will prove that there is no basis for the
uncountable regular Hausdorff spaces of cardinality $\aleph_1$.
Each of these results gives a strong refutation of a
well known and longstanding conjecture.
The results all stem from the analysis of oscillations of
coherent sequences $\<e_\alpha:\alpha < \omega_1\>$
of finite-to-one functions.
I expect that the methods presented will have other applications as well.
\end{abstract}

\maketitle

\section{Introductions}
One goal in mathematics is to classify structures, particularly those
which lie at the core of mathematics ---
graphs, groups, orders, relations, solids,
surfaces, vector spaces, and topologies.\footnote{All topological
spaces in this paper are regular and Hausdorff.}
Classifications can come in many forms.
The strongest and most infrequent is a complete and informative
listing --- up to presentation ---
of the structures of a given class of objects.
The following are examples of such a classification.
\begin{thm}
There are five regular polyhedra:
the tetrahedron, the cube, the octahedron, the dodecahedron,
and the icosahedron.
\end{thm}
\begin{thm}
Any two closed orientable 2-manifolds with the same genus are
homeomorphic.
\end{thm}
\begin{thm} \cite{oriented_systems}
Any countable transitive relation is Tukey equivalent to
one of the following forms
for some non-negative integers $n_i$ $(i < 3)$:
$$n_0 \cdot 1 \oplus n_1 \cdot \omega \oplus n_2 \cdot \seq.$$
\end{thm}

Frequently it is not possible to have a reasonable
complete classification of this sort.
If this happens, it is still possible to ask whether the objects can
be classified from ``above'' or from ``below.''
In each case, there is a need for a notion of reduction $\leq$
(e.g. being a minor, a normal subgroup, subgraph, a subspace, etc.).
A family $\Bcal$ is a \emph{basis} for a class $\Ccal$
if for any member $C$ of $\Ccal$, there is a $B$ in $\Bcal$ such that
$B \leq C$.
A family $\Ucal$ is a \emph{universal} for $\Ccal$ if any member $C$ of
$\Ccal$, there is a $U$ in $\Ucal$ such that $C \leq U$.
An example of a basis theorem is the following result of Kuratowski.
\begin{thm} \cite{nonplanar_basis}
If $G$ is a non-planar graph, then $G$ contains either $K_{3,3}$ or
$K_5$ as a minor.
\end{thm}
The classification of finite simple groups can also be considered
a basis theorem for the class of finite groups with reduction
being ``is contained as a normal subgroup.''
The following are two examples of universality results.
\begin{thm}
$\Rbb$ is a universal separable linear order.
\end{thm}
\begin{thm}
The irrationals are surjectively universal for the
separable completely metrizable spaces.\footnote{Unfortunately
it is not even well known in many mathematical circles
that the irrationals are themselves completely metrizable.}
\end{thm}

In this paper we will focus on basis problems.
Since the whole class $\Ccal$ is trivially a basis,
the goal is to find a basis which is informative.
This is often quantified in terms of restrictions on
its cardinality.
In cases where the basis must be infinite, it is also of
interest to know whether it can contain only minimal
elements, hence avoiding any redundancies.

In some cases it is possible to show that there is no
reasonable basis --- that any such basis is either large in
terms of cardinality or must contain redundancies.
While such a result is not exactly in line with the goal of
classification, it is still informative as to the limits
of our ability to understand the class. 

A phenomenon of this sort has been repeatedly observed for objects
at the level of the continuum.
The following classical result of Sierpi\'nski is a typical example.
\begin{thm} \cite{many_iso_types}
There is a subset $X \subseteq \Rbb$ of cardinality continuum
such that if $f \subseteq X^2$ is a continuous injection, then
$f$ differs from the identity map on a set of size less than continuum.
In particular, the suborders of $X$ witness that
any basis for the separable linear orders of
cardinality continuum must be uncountable and contain redundancies.
\end{thm}

In this paper we will focus on basis problems for uncountable objects.
For finite and countable objects, classification theorems
are typically decidable by the accepted axioms of mathematics ($\ZFC$).
This phenomenon is largely explained by Shoenfield's absoluteness
theorem \cite{Shoenfield_absoluteness}.
When dealing with uncountable objects, however, it is often the case
that additional axioms are necessary if any classification is
consistent with $\ZFC$.
This is because the non-existence of reductions
between fixed mathematical objects can be fail to be absolute.
Axioms such as Jensen's $\diamond$ and the
Continuum Hypothesis allow one to inductively define examples while
diagonalizing against a list of all possible reductions.
If $2^{\aleph_0} = \aleph_1$, then the pathologies mentioned above
occur already at the first uncountable cardinal.
For instance, Sierpi\'nski's theorem easily yields
the following fact.
\begin{thm}
If the Continuum Hypothesis is true, then
any basis for the uncountable linear or the uncountable topological
spaces must have cardinality at least $2^{\aleph_1}$ and contain
redundancies.
\end{thm}
On the other hand, the lack of embeddings between two uncountable
sets of reals is not absolute as the following result of Baumgartner's
demonstrates.
\begin{thm}\cite{reals_isomorphic}
If $X$ and $Y$ are two uncountable sets of reals, then there is
a forcing which preserves uncountability and which
adds an order preserving map from $X$ into $Y$.
Moreover, the Proper Forcing Axiom implies that
such an embedding already exists provided $X$ has cardinality
$\aleph_1$.
\end{thm}
The Proper Forcing Axiom (PFA) is a statement which allows one
to find embeddings \emph{generically} between two objects of cardinality
$\aleph_1$.
The use of this axiom can be likened to that of Erd\"os's
probabilistic method in finite combinatorics
(see \cite{probabilistic_method}).
The main difference is that the notion of a probability space is
abstracted and an axiom is required for
the deduction
``if an object can be selected
with positive probability, then it exists.''
This axiom is natural to assume if one wishes to classify objects at
the level of $\aleph_1$, the smallest uncountable cardinal.
In fact, there is the following empirical fact:\footnote{As
with any claim of this sort, there are counterexamples.
A notable example is Kat\v{e}tov's problem --- see \cite{katprob}.}
For a given class $\Ccal$ consisting of
structures of size $\aleph_1$ the existence of a classification or
a basis is either a consequence
of $\PFA$ or else is provably false without additional axiomatic
assumptions.
This has been partially explained by the work of Woodin \cite{Pmax}.

In his 1998 address to the International Congress of Mathematicians,
Todor\v{c}evi\'{c} outlined a number of basis conjectures\footnote{
Todor\v{c}evi\'{c} lists these as conjectures in \cite{basis_problems}.
This is not to say that the author of the statement necessarily
speculated a certain resolution.}
--- in the presence of $\PFA$ ---
which had been under consideration since the 1960's and 70's.

\begin{conjecture}\label{Tukey_conjecture}
(Todor\v{c}evi\'{c}; \cite{basis_problems})
If $R$ is a binary relation, then either
$R \leq \aleph_0 \cdot \omega_1$ or
$[\omega_1]^{< \aleph_0} \leq R$.
\end{conjecture}

\begin{conjecture}\label{SL_conjecture}
(Hajnal, Juhasz; \cite{alphaS_alphaL})
If $X$ is a regular Hausdorff space, then the following are
equivalent:
\begin{enumerate}

\item $X$ is hereditarily separable.

\item $X$ is hereditarily Lindel\"of.

\item $X$ does not contain an uncountable discrete subspace.

\end{enumerate}
\end{conjecture}

\begin{conjecture} \label{rectangle}
(Galvin\footnote{This attribution is a conjecture of Hajnal.
The first time this appeared in print seems to be
\cite{embedding_neg_relation}.})
The partition relation $\omega_1 \rightarrow (\omega_1;\omega_1)^2_2$
holds.
\end{conjecture}

\begin{conjecture} \label{Shelah_conjecture}
(Shelah; \cite{Countryman:Shelah})
The uncountable linear orders have a five element basis consisting
of $X$, $\omega_1$, $\omega_1^*$, $C$, and $C^*$
where $X$ is a set of reals of cardinality $\aleph_1$ and $C$ is a
Countryman type.\footnote{A linear order
$C$ is Countryman if it is uncountable
and yet $C^2$ is the union of countably many chains in the
coordinate-wise order.}
\end{conjecture}

\begin{conjecture}\label{Gruenhage_conjecture}
(Gruenhage; \cite{open_problems_topology:Gruenhage})
The uncountable regular Hausdorff spaces have a three element basis
consisting of a set of reals of cardinality $\aleph_1$ with
the metric, the Sorgenfrey, and the discrete topology.
\end{conjecture}

\begin{conjecture}\label{Fremlin_problem}
(Fremlin; DN of \cite{problems:Fremlin})
If $X$ is a compact space in which closed sets
are $G_\delta$, then there is an at most two-to-one
map from $X$ into a compact metric space.
\end{conjecture}

Todor\v{c}evi\'{c} has, assuming $\PFA$,
proved that Conjecture \ref{Tukey_conjecture}
for the transitive relations \cite{class_trans},
and confirmed
the equivalence of (1) and (3) in Hajnal and Juhasz's problem
\cite{forcing_partition}.
In \cite{linear_basis}, I confirmed that Conjecture \ref{Shelah_conjecture}
follows from $\PFA$.
Gruenhage has proved that Conjecture \ref{Gruenhage_conjecture}
holds for the class of cometrizable spaces assuming $\PFA$
\cite{cosmicity}.
Todor\v{c}evi\'c has verified that Conjecture \ref{Fremlin_problem} holds
for the Rosenthal compacta in $\ZFC$ \cite{Rosenthal_cpt}.

In this paper I will prove results in the other direction, refuting
all of the conjectures except
Conjectures \ref{Shelah_conjecture} and \ref{Fremlin_problem}.
\begin{thm}\label{no_rectangle}
There is a function $f:[\omega_1]^2 \to \omega_1$ such that
whenever $A,B \subseteq \omega_1$ are uncountable and
$\xi < \omega_1$ there are $\alpha$ in $A$ and $\beta$ in
$B$ such that $\alpha < \beta$ and $f (\alpha,\beta) = \xi$.
In particular $\omega_1 \not \rightarrow (\omega_1;\omega_1)^2_2$.
\end{thm}

\begin{thm} \label{Tukey_thm}
There is a family $\Fcal$ of binary relations
which contains an antichain of cardinality $2^{\aleph_1}$, is downwards
directed and such
that for all $R$ in $\Fcal$,
$R \not \leq \aleph_0 \cdot \omega_1$ and
$[\omega_1]^{<\aleph_0} \not \leq R$.
\end{thm}

\begin{thm}\label{Lspace}
There is a non-separable, hereditarily Lindel\"of topological space.
\end{thm}

\begin{thm} \label{no_top_basis}
Any basis for the uncountable topological
spaces has cardinality strictly greater than $\aleph_1$.
\end{thm}

It is worth noting that the existence of an L space alone gives a
refutation of these conjectures.
The above results, however, take the failures of these conjectures
a step further and require slightly different combinatorial refinements.

The theorems in this paper are a consequences of the following result
and its technical strengthenings.
The reader is referred to Section \ref{trace_section} for the
undefined notions.
\begin{thm}
Suppose that $\<e_\beta:\beta < \omega_1\>$ is a coherent sequence
of finite-to-one functions with $e_\beta:\beta \to \omega$ for each
$\beta < \omega_1$ and suppose that $L:[\omega_1]^2 \to [\omega_1]^{<\aleph_0}$
satisfies the properties of a lower trace function.
If $A,B \subseteq \omega_1$ are uncountable, then the set of
integers of the form
$$\osc(e_\alpha \restriction L(\alpha,\beta),
e_\beta \restriction L(\alpha,\beta))$$
for $\alpha < \beta$ with $\alpha$ in $A$, $\beta$ in $B$ contains
arbitrarily long intervals.
\end{thm}

This can be likened to the following result of Todor\v{c}evi\'c, which he
used to draw a number of conclusions about Conjecture \ref{SL_conjecture}.
\begin{thm}\cite{partition_problems}
If $X$ and $Y$ are unbounded and countably directed in
$(\ww,<^*)$ and consist of monotonic functions, then there
is an natural number $l$ such that for all $n$ there are 
$x$ and $y$ in $X$ and $Y$ respectively such that
$$\osc(x,y) = n + l.$$
\end{thm}

The paper is organized as follows.
Section \ref{trace_section} provides some background on the method
of minimal walks and introduces the lower trace function which is
used in the statement of the main technical theorem on oscillations.
Section \ref{models_section} provides a review of elementary
submodels.
The main technical results of the paper are proved
in Section \ref{osc_section}.
Theorem \ref{no_rectangle} is then deduced in Section \ref{o_section}
and a coloring $o:[\omega_1]^2 \to \omega$ is introduced which
harnesses most of the strength of the more technical main theorem.
Theorem \ref{Tukey_thm} is then deduced in Section \ref{Tukey_section}
after some motivation is provided.
The paper closes with Section \ref{Lspace_section} where
Theorems \ref{Lspace} and \ref{no_top_basis} are proved.
This final section also contains a basic analysis of the L space
which is constructed.

This paper is intended to be accessible to any interested reader
who is fluent in set theory --- basic background can be found in,
e.g., \cite{set_theory:Jech}, \cite{set_theory:Kunen}.
Elementary submodels will be employed in a number of points in the
argument.
Ironically this represents the only non-elementary technique used
in the proof.
The essentials are reviewed in Section \ref{models_section}.
The reader is referred to III.1 of \cite{multiple_forcing} for more
information on elementary submodels and stationary sets.
The proofs will also employ the method of minimal walks introduced in
\cite{partitioning_ordinals}.
The necessary background is presented in Section \ref{trace_section}.
The reader is referred to \cite{cseq} for further reading on minimal
walks.

The notation is fairly standard.
All ordinals are von Neumann ordinals --- they are the set
consisting of their
predecessors.
In particular, $n = \{0,1,\ldots,n-1\}$ and 
the first infinite ordinal $\omega$
is the set of all finite ordinals $\{0,1,2,\ldots\}$ and is identified with
the natural numbers.
All counting starts at 0.
If $k$ is a natural number and $X$ is a set, then $[X]^k$ is the set of all
$k$-element subsets of $X$.
If $X$ has a canonical linear ordering associated with it and $a$ is
in $[X]^k$, then $a$ will be identified with the increasing sequence
which enumerates it.
If $a$ and $b$ are finite subsets of $\omega_1$, then
$a < b$ will be used to abbreviate ``$\alpha < \beta$ whenever $\alpha$
is in $a$ and $\beta$ is in $b$.''
Similarly one defines statements such as $\alpha < b$ and
$a < \beta$ if $\alpha$ and $\beta$ are ordinals.

\section{The trace functions}
\label{trace_section}
In this section I will provide the necessary background on minimal
walks.
Minimal walks are facilitated by a \emph{$C$-sequence} or ladder
system which one uses to ``walk'' from an ordinal $\beta$
down to a smaller ordinal $\alpha$.
Here a $C$-sequence is a sequence $\<C_\alpha:\alpha < \omega_1\>$
such that
$C_\alpha$ is a cofinal subset of $\alpha$ and if $\gamma < \alpha$
then $C_\alpha \cap \gamma$ is finite.
It will be useful at certain points to assume that 0 is an element of
every $C_\alpha$.

The following two functions will be of interest to us.
The upper trace will not be necessary, but is useful in making
the other definitions more transparent.

\begin{defn} \cite{cseq} (upper trace)
If $\alpha < \beta$, then $\tr(\alpha,\beta)$ is defined
recursively by
$$\tr(\alpha,\alpha) = \emptyset,$$
$$\tr(\alpha,\beta) =
\tr\typezero(\alpha,\min(C_\beta \setminus \alpha)\typezero) \cup \{\beta\}.$$
\end{defn}
\begin{defn} (lower trace)
If $\alpha < \beta$, then $L(\alpha,\beta)$ is defined
recursively by
$$L(\alpha,\alpha) = \emptyset,$$
$$L(\alpha,\beta) =
L\typezero(\alpha,\min(C_\beta \setminus \alpha)\typezero) \cup \{\max(C_\beta \cap
\alpha)\}
\setminus \max(C_\beta \cap
\alpha).$$
\end{defn}
Hence the upper trace $\tr(\alpha,\beta)$ is enumerated by the sequence
$\beta = \beta_0 > \beta_1 > \cdots > \beta_{l-1}$
where $\beta_{i+1} = \min(C_{\beta_i} \setminus \alpha)$ and
$\alpha = \beta_l$ is in $C_{\beta_{l-1}}$.
It is easily verified that
the lower trace $L(\alpha,\beta)$ is listed as
$\xi_0 \leq \xi_1 \leq \cdots \leq \xi_{l-1}$
where $\xi_i$ is the maximum of
$$\bigcup_{j=0}^{i} C_{\beta_j} \cap \alpha.$$
While the upper trace is well studied, this is, to my knowledge, the
first analysis of the lower trace (the \emph{full lower trace}
$F(\alpha,\beta)$
is a different matter --- see \cite{cseq}).

The lower trace can be axiomatized by
the following facts.

\begin{fact} \label{L_trans}
If $\alpha \leq \beta \leq \gamma$ and $L(\beta,\gamma) < L(\alpha,\beta)$,
then
$$L(\alpha,\gamma) = L(\alpha,\beta) \cup L(\beta,\gamma).$$
\end{fact}

\begin{proof}
Let $\alpha$, $\beta$, and $\gamma$ satisfy the hypothesis of the
fact.
Observe that $L(\beta,\gamma) < \alpha$ and hence
$C_\zeta \cap \alpha = C_\zeta \cap \beta$ whenever
$\zeta$ is in $\tr(\beta,\gamma)$.
Hence $\beta$ is in $\tr(\alpha,\gamma)$ and
$\tr(\alpha,\gamma) = \tr(\alpha,\beta) \cup \tr(\beta,\gamma)$.
Let $$\gamma = \gamma_0 > \gamma_1 > \ldots > \gamma_{l -1}$$
enumerate $\tr(\alpha,\gamma)$ and let $l_0$ be such that
$\beta = \gamma_{l_0}$.
Then $L(\alpha,\gamma)$ is listed as
$$\xi_0 \leq \xi_1 \leq \cdots \leq \xi_{l-1}$$
where $$\xi_i = \max \bigcup_{j=0}^i C_{\gamma_j} \cap \alpha.$$
If $j < l_0$, then $C_{\gamma_j} \cap \alpha = C_{\gamma_j} \cap
\beta$
and so
$L(\beta,\gamma) = \{\xi_j\}_{j=0}^{l-1}$.
On the other hand,
$$\max\typefour(C_{\gamma_{l_0}} \cap \alpha\typefour) > \xi_{l_0-1}$$
and so
$$\max \bigcup_{j=0}^i \typefour( C_{\gamma_j} \cap \alpha \typefour) = 
\max \bigcup_{j=l_0}^i \typefour( C_{\gamma_j} \cap \alpha \typefour)$$
and hence $L(\alpha,\beta) = \{\xi_j\}_{j=l_0}^{l-1}$.
\end{proof}

\begin{fact} \label{L_lim}
If $\delta$ is a limit ordinal, then
$$\lim_{\xi \to \delta} \min L(\xi,\delta) = \delta.$$
\end{fact}

\begin{proof}
This follows immediately from the
observation that
$$\min L(\xi,\delta) = \max(C_\delta \cap \xi).$$
\end{proof}

For the stronger coloring $o$ used later in the paper, we will also
need a fixed enumeration $\<w_\xi:\xi < \omega_1\>$ of
$C(2^\omega,\omega)$ which lists all elements stationarily often.
Here $C(2^\omega,\omega)$ is the set of all continuous functions
from $2^\omega$ into $\omega$.
Since $2^\omega$ is compact and $\omega$ is discrete, this collection
is countable and its elements can be thought of as weighted clopen
partitions of $2^\omega$.
We will also need a fixed sequence $z_\alpha$ $(\alpha < \omega_1)$
of distinct elements of $2^\omega$.

This allows for the following definition.
\begin{defn}
If $\alpha < \beta$, then
$\mu(\alpha,\beta)$ is the function defined on
$L(\alpha,\beta)$ recursively by
$$\mu\typetwo(\alpha,\beta\typetwo)\typezero(\max(C_\beta \cap
\alpha)\typezero) = w_\beta,$$
$$\mu(\alpha,\beta)(\gamma) =
\mu\typezero(\alpha,\min(C_\beta \setminus \alpha)\typezero)(\gamma)$$
for all $\gamma > \max(C_\beta \cap \alpha)$ in $L(\alpha,\beta)$.
Let $\mu(\alpha,\beta;\xi)$ be the pointwise evaluation of $\mu(\alpha,\beta)$
at $z_\xi$.
\end{defn}
It can be checked that $\mu(\alpha,\beta)$ is the mapping
defined by $\xi_{i} \mapsto w_{\beta_i}$ if $\xi_{i-1} < \xi_{i}$ or
$i = 0$ where $\xi_i$ and $\beta_i$ are as above.
Like $L$, the function $\mu$ can be axiomatized by the following
properties.
\begin{fact}\label{mu_trans}
If $\alpha < \beta < \gamma$ and $L(\beta,\gamma) < L(\alpha,\beta)$,
then
$$\mu(\alpha,\gamma) = \mu(\alpha,\beta) \cup \mu(\beta,\gamma)$$
(i.e. $\mu(\alpha,\beta)$ and $\mu(\beta,\gamma)$ are both
restriction of $\mu(\alpha,\beta)$.).
\end{fact}
\begin{fact}
\label{mu_lim}
If $\xi < \delta$, then
$\mu\typetwo(\xi,\delta\typetwo)
\typezero(\min L(\xi,\delta)\typezero) = w_\delta$.
\end{fact}

The main technical result of the paper will be concerned with
counting oscillations between elements of a coherent\footnote{A family
of functions is \emph{coherent} in this context if
any two differ on a finite set on the intersection of their domains.}
sequence $\<e_\beta:\beta < \omega_1\>$ where $e_\beta$ is a
finite-to-one function from $\beta$ into $\omega$.
The following construction gives a standard example of a coherent
family of finite-to-one functions.
\begin{defn} \cite{partitioning_ordinals}(maximal weight)
If $\alpha \leq \beta$, then 
$\rho_1(\alpha,\beta)$ is defined recursively by
$$\rho_1(\alpha,\alpha) = 0$$
$$\rho_1(\alpha,\beta) = \max\Big(|C_\beta \cap \alpha|,
\rho_1\typezero(\alpha,\min(C_\beta \setminus \alpha)\typezero)\Big).$$
\end{defn}
Alternately, $\rho_1(\alpha,\beta)$ is the maximum value
of the form $|C_\xi \cap \alpha|$ where $\xi$ ranges over
$\tr(\alpha,\beta)$.
Put $e_\beta(\alpha) = \rho_1(\alpha,\beta)$.
\begin{fact}\cite{partitioning_ordinals}\label{rho1_fact}
$\<e_\beta:\beta < \omega_1\>$ is a coherent sequence of finite-to-one
functions.
\end{fact}
\begin{proof}
Let $\beta \leq \beta' < \omega_1$ and $n < \omega$ be given and
set $D$ equal to the set of all $\alpha < \beta$ such that either
$e_\beta(\alpha) \leq n$ or $e_\beta(\alpha) \ne
e_{\beta'}(\alpha)$.
It suffices to show that $D$ has no limit points.

To this end,
suppose $\delta \leq \beta$.
It is easy to check that there is a $\delta_0 < \delta$ such that
$$\tr(\alpha,\beta) = \tr(\alpha,\delta) \cup \tr(\delta,\beta),$$
$$\tr(\alpha,\beta') = \tr(\alpha,\delta) \cup \tr(\delta,\beta'),$$
$$|C_\delta \cap \alpha| > n$$ 
whenever $\delta_0 < \alpha < \delta$.
If $\delta_0 < \alpha < \delta$, then
$$e_{\beta}(\alpha) \geq |C_\delta \cap \alpha| > n$$
$$e_{\beta}(\alpha) = e_\delta(\alpha) = e_{\beta'}(\delta)$$
and hence $\alpha$ is not in $D$.
Consequently $\delta$ is not a limit point of $D$.
\end{proof}
\begin{remark}
It is interesting to note that while coherence and the finite-to-one
property are at tension with each other, the verifications of these
properties in the previous fact are virtually identical.
\end{remark}
Throughout this paper we will now assume that $L$, $\mu$,
$\<e_\beta:\beta < \omega_1\>$,
$\<w_\alpha:\alpha <\omega_1\>$,
and $\<z_\alpha:\alpha < \omega_1\>$
are fixed objects defined as above.

\section{Basic facts about elementary submodels}
\label{models_section}
In the proof of the main theorem,
we will need the following facts about elementary submodels
of $H(\aleph_2)$.
Recall that $H(\aleph_2)$ is the collection of all sets of hereditary
cardinality less than $\aleph_2$.
This structure is of interest to us since it contains all of
the objects relevant to us
and satisfies all of the axioms of $\ZFC$ except the power set axiom.
An elementary submodel $M$ of $H(\aleph_2)$ is a subset of
$H(\aleph_2)$ such that whenever $\phi$ is a logical formula with parameters
in $M$, $M$ satisfies $\phi$ if and only if $H(\aleph_2)$ satisfies
$\phi$.
The following fact is well known.

\begin{fact} \label{omega1_fact}
If $M$ is a countable elementary submodel of $H(\aleph_2)$, then
$M \cap \omega_1$ is an ordinal.
Moreover, if
$F$ is a finite subset of $H(\aleph_2)$, then the set
of all ordinals of the form $M \cap \omega_1$ such that
$M$ is a countable elementary submodel of $H(\aleph_2)$ contains
a club.\footnote{A subset of $\omega_1$ is \emph{club} if
it is closed and unbounded.
A subset of $\omega_1$ is \emph{stationary} if it intersects every
club.}
\end{fact}

The following standard facts are very useful when working with
countable elementary submodels.
\begin{fact} \label{uncountable_test1}
If $M$ is a countable elementary submodel of $H(\aleph_2)$ which
contains some element $X$, then $X$ is countable
iff $X \subseteq M$.
\end{fact}

\begin{fact} \label{uncountable_test2}
If $M$ is an elementary submodel of $H(\aleph_2)$ which contains as an element
some subset $A$ of $\omega_1$, then $A$ is uncountable iff
$A \cap M \cap \omega_1$ is unbounded in $M \cap \omega_1$.
\end{fact}

\begin{fact} \label{definable_is_in}
If $M$ is an elementary submodel of $H(\aleph_2)$, $X$ is in $H(\aleph_2)$
and $X$ is definable from a formula with parameters in $M$, then
$X$ is in $M$.
\end{fact}

\section{Oscillations on the lower trace}
\label{osc_section}
In this section I will prove main technical lemmas of the paper.
We will take the following as our definition of the oscillation map.
\begin{defn}
Suppose that $s$ and $t$ are two functions defined on
a common finite set of ordinals $F$.
Let $\Osc(s,t;F)$ be the set of all $\xi$ in $F$ such that
$s(\xi^-) \leq t(\xi^-)$ and $s(\xi) > t(\xi)$ where
$\xi^-$ is the greatest element of $F$ less than $\xi$.\footnote{Actually
we will later see that we are really interested in counting
oscillations between the relations $=$ and $\ne$
rather than the more conventional $<$ and $>$.
The above definition is a compromise between the two.}
\end{defn}

The following notation will be convenient.
 
\begin{defn}
If $\alpha < \beta < \omega_1$,
let $\Osc(\alpha,\beta)$ denote
$$\Osc\typeone(e_\alpha,e_\beta;L(\alpha,\beta)\typeone)$$
and $\osc(\alpha,\beta)$ denote the cardinality of $\Osc(\alpha,\beta)$.
\end{defn}

\begin{lem} \label{osc_lem}
For every uncountable pairwise disjoint $\Acal \subseteq [\omega_1]^k$
and $\Bcal \subseteq [\omega_1]^l$ and all $w$ in
$C(2^\omega,\omega)$, there are
$b_m$ $(m < \omega)$ in $\Bcal$
such that for all $n$ there is an $a$ in $\Acal$ and
$\{\xi_i:i < n\}$ such that
for all $m \leq n$, $i < k$, $j < l$ the following hold:
\begin{enumerate}

\item  $a < b_m$.

\item $\Osc\typeone(a(i),b_m(j)\typeone)$ is the disjoint union of
$\Osc\typeone(a(i),b_0(j)\typeone)$ and $\{\xi_{m'}:m' < m\}$.

\item
$\mu\typeone(a(i),b_0(j)\typeone)$ is $\mu\typeone(a(i),b_m(j)\typeone)$
restricted to $L\typeone(a(i),b_0(j)\typeone)$.

\item
$\mu\typeone(a(i),b_m(j)\typeone)(\xi_{m'}) = w$ whenever $m' < m$.

\end{enumerate}
\end{lem}

\begin{remark}
While the generality of this theorem will be needed to derive
some of the results in the paper, many of them can get by with
weaker formulations.
For instance, the analysis of $\osc$
only requires the theorem for $k = l = 1$
and only conclusions (1) and
(2).
This is sufficient to prove Theorem \ref{no_rectangle}
The later conditions are needed for the properties of the
stronger coloring $o$ presented in Section \ref{o_section}.
\end{remark}
Lemma \ref{osc_lem} will be derived from the following theorem.

\begin{lem} \label{osc_recursion}
Let $\Acal \subseteq [\omega_1]^k$ and $\Bcal \subseteq [\omega_1]^l$
be uncountable and pairwise disjoint.
There is a closed an unbounded set
of $\delta < \omega_1$ such that
if $a$ is in $\Acal \setminus \delta$,\footnote{It will be convenient
to let $\Acal \setminus \delta$ denote the set of all $a$ in 
$\Acal$ such that $\delta \leq a$ whenever $\Acal$ is a collection
of finite subsets of $\omega_1$ and $\delta < \omega_1$.}
$b$ is in $\Bcal \setminus \delta$,
and $R$ is in $\{=,>\}$,
then there are $a^+$ in $\Acal \setminus \delta$,
and $b^+$ in $\Bcal \setminus \delta$ such that for all $i < k$ and
$j < l$:
\begin{enumerate}

\item
The maximum of
$L\typeone(\delta,b(j)\typeone)$
is less than both $\Delta\typefour(e_{a(i)},e_{a^+(i)}\typefour)$ and 
$\Delta\typefour(e_{b(j)},e_{b^+(j)}\typefour)$.

\item There is a non-empty $L^+$ which does not depend on $j$ such that
$L\typeone(\delta,b(j)\typeone) < L^+$ and
$L\typeone(\delta,b^+(j)\typeone) = L\typeone(\delta,b(j)\typeone) \cup L^+$.

\item
If $\xi$ is in $L^+$, then
$e_{a^+(i)}(\xi)\ R\ e_{b^+(j)}(\xi)$.

\item
$\mu\typeone(\delta,b(j)\typeone)$ is $\mu\typeone(\delta,b^+(j)\typeone)$ restricted to
$L\typeone(\delta,b(j)\typeone)$.

\item
$\mu\typeone(\delta,b^+(j)\typeone)(\min L^+) = w_\delta$.
 
\end{enumerate}
\end{lem}

\begin{proof}
Let $M$ be an elementary submodel of $H(\aleph_2)$ which contains
everything relevant and let $\delta = M \cap \omega_1$.
By Fact \ref{omega1_fact},
it suffices to show that $\delta$ satisfies
the conclusion of the theorem.

First suppose that $R$ is $=$.
Applying Fact \ref{rho1_fact}, find a
$\gamma_0 < \delta$ satisfying the following conditions:
\begin{enumerate}

\item If $j < l$, then $L\typeone(\delta,b_m(j)\typeone) < \gamma_0$.

\item If $i < k$, $j < l$, and $\gamma_0 < \xi < \delta$, then
$e_{a(i)}(\xi) = e_{b(j)}(\xi)$.

\end{enumerate}
Applying Facts \ref{L_lim} and \ref{mu_lim},
pick a $\gamma < \delta$ such that if $\gamma < \xi < \delta$, then
$\gamma_0 < L(\xi,\delta)$.
Consider the set $D \subseteq \omega_1$
consisting of all $\delta^+$ such that for some
$a^+$ in $\Acal \setminus \delta^+$
and $b^+$ in $\Bcal \setminus \delta^+$ the following conditions
are satisfied:
\begin{enumerate}

\item
$e_{a^+(i)} \restriction \gamma_0 = e_{a(i)} \restriction \gamma_0$
and
$e_{b^+(j)} \restriction \gamma_0 =
e_{b(j)} \restriction \gamma_0$
for all $i < k$ and $j < l$.

\item If $\gamma < \xi < \delta^+$, then
$\gamma_0 < L(\xi,\delta^+)$ and
$\mu\typetwo(\xi,\delta^+\typeone)\typeone(\min
L(\xi,\delta^+)\typeone)
= w_{\delta^+} = w_{\delta}$.

\item If $\gamma_0 < \xi < \delta^+$, $i < k$, and $j < l$, then
$e_{a^+(i)}(\xi) = e_{b^+(j)}(\xi)$.

\item $\mu\typeone(\delta^+,b^+(j)\typeone) =
\mu\typeone(\delta,b(j)\typeone)$.

\end{enumerate}
Observe that for all $\beta \geq \gamma_0$,
$e_{\beta} \restriction \gamma_0$ is in $M$
since by Fact \ref{rho1_fact} it
differs from $e_{\gamma_0}$ on a finite set.
Hence $D$ is definable from the parameters
$\gamma_0$, $\gamma$, $e_{a(i)} \restriction \gamma_0$, etc. which are
all elements of $M$.
Therefore $D$ is in $M$ by Fact \ref{definable_is_in}.
Since $D$ has $\delta$ as a member,
it is uncountable by Fact \ref{uncountable_test1}.
Hence there is a $\delta^+ > \delta$ in $D$.
Let $a^+ \in \Acal \setminus \delta^+$ and $b^+ \in \Bcal \setminus \delta^+$
witness that $\delta^+$ is in $D$.

Now let $i < k$ and $j < l$ be arbitrary.
First observe that
$$\gamma_0 \leq \Delta\typefour(e_{a(i)},e_{a^+(i)}\typefour),$$
$$\gamma_0 \leq \Delta\typefour(e_{b(j)},e_{b^+(j)}\typefour).$$
Put $L^+ = L(\delta,\delta^+)$.
Notice that
$L\typeone(\delta,b(j)\typeone) < L^+$ and hence
$$L\typeone(\delta,b^+(j)\typeone) =
L\typeone(\delta^+,b^+(j)\typeone) \cup L^+$$
holds by Fact \ref{L_trans}.
Since $\delta < \delta^+$, $L^+$ is non-empty.
If $\xi$ is in $L^+$,
then $\gamma_0 < \xi < \delta^+$ and so
$$e_{a^+(i)}(\xi)\ R\ e_{b^+(j)}(\xi).$$

Now suppose that $R$ is $>$.
Put $E$ equal to the set of all limit $\nu < \omega_1$ such that
for all $a_0$ in $\Acal \setminus \nu$, all
$\nu_0 < \nu$, $\epsilon < \omega_1$, $n < \omega$, and
finite $L^+ \subseteq \omega_1 \setminus \nu$
there is an
$a_1$ in $\Acal \setminus \epsilon$ with
$$\nu_0 \leq \Delta\typefour(e_{a_0(i)},e_{a_1(i)}\typefour),$$
$$e_{a_1(i)}(\xi) > n$$
whenever $i < k$ and $\xi$ is in $L^+$.
By Fact \ref{definable_is_in}, $E$ is an element of $M$.
\begin{claim} \label{E_claim}
$\delta$ is in $E$.
In particular, $E$ is uncountable by Fact \ref{uncountable_test1}.
\end{claim}

\begin{proof}
Let $a_0$, $\nu_0$, $\epsilon$, $n$, and $L^+$ be given as in the definition
of $E$ for $\nu = \delta$.
By Fact \ref{rho1_fact} we may assume without loss of generality
that $\nu_0$ is an upper bound
for all $\xi < \delta$ such that $e_{a_0(i)}(\xi) \leq n$ for some $i < k$.
Now, applying elementarity of $M$, there is a
$\delta^+$ above $\epsilon$, $\delta$, and $\max L^+$
and an $a_1$ in $\Acal \setminus \delta^+$ such that the following conditions
are satisfied:
\begin{enumerate}

\item For all $i < k$,
$e_{a_0(i)} \restriction \nu_0 = e_{a_1(i)} \restriction \nu_0$.

\item If $\nu_0 < \xi < \delta^+$ and $i < k$, then
$e_{a_1(i)}(\xi) > n$.

\end{enumerate}
Since $L^+ \subseteq \delta^+ \setminus \delta$,
this completes the proof of the claim.
\end{proof}
Applying elementarity of $M$ and uncountability of $E$,
find an element $\gamma_0$ of $E$ such that
$L\typeone(\delta,b(j)\typeone) < \gamma_0 < \delta$ for all $j < l$.
Applying Fact \ref{L_lim},
find a $\gamma < \delta$ such that if $\gamma < \xi < \delta$,
then
$\gamma_0 < L(\xi,\delta)$.
Again using elementarity of $M$, select
a limit $\delta^+ > \delta$ and
a $b^+$ in $\Bcal \setminus \delta^+$ so that the following conditions
are satisfied:
\begin{enumerate}

\item $e_{b^+(j)} \restriction \gamma_0 =
e_{b(j)} \restriction \gamma_0$ for all $j < l$.

\item If $\gamma < \xi < \delta^+$, then
$\gamma_0 < L(\xi,\delta^+)$ and
$\mu\typetwo(\xi,\delta^+\typeone)\typeone(\min
L(\xi,\delta^+)\typeone)
= w_{\delta^+} = w_{\delta}$.

\item $\mu\typeone(\delta^+,b^+(j)\typeone) =
\mu\typeone(\delta,b(j)\typeone)$.

\end{enumerate}
Put $L^+ = L(\delta,\delta^+)$.
Applying the definition of $E$, find an $a^+$ in
$\Acal \setminus \delta$ such that for all $i < k$, $j < l$, and $\xi$
in $L^+$
$$L\typeone(\delta,b(j)\typeone) < \Delta\typefour(e_{a(i)},e_{a^+(i)}\typefour)$$
$$e_{a^+(i)}(\xi) > e_{b^+(j)}(\xi).$$
The rest of the verification is as in the previous case.
This completes the proof of Lemma \ref{osc_recursion}.
\end{proof}

Now we are ready to prove Lemma \ref{osc_lem}
\begin{proof} 
Let $\Acal$ and $\Bcal$ be given and select a countable elementary
submodel $M$ of $H(\aleph_2)$ containing everything relevant such that
$w_\delta = w$ where $\delta = M \cap \omega_1$.
Since $M$ contains $\Acal$ and $\Bcal$, the closed and unbounded
set provided by Lemma \ref{osc_recursion} is in $M$ and therefore
$\delta$ is in this closed an unbounded set.
Using Lemma \ref{osc_recursion}, it is easy to select
$a_m$ $(m < \omega)$ in $\Acal \setminus \delta$ and
$b_m$ $(m < \omega)$ in $\Bcal \setminus \delta$
so that for all $m < \omega$ the following conditions are satisfied whenever
$i < k$ and $j < l$:
\begin{enumerate}

\item $L\typeone(\delta,b_{m+1}(j)\typeone)$ is a proper initial part of
$L\typeone(\delta,b_m(j)\typeone)$.

\item The minimum of $L\typeone(\delta,b_{m+1}(j)\typeone) \setminus L\typeone(\delta,b_m(j)\typeone)$
is an ordinal $\xi_m$ which does not depend on $j$.

\item $\Osc\typeone(a_{m+1}(i),b_{m+1}(j);L(\delta,b_{m+1}(j))\typethree)$
is formed by adding $\xi_m$ to 
$\Osc\typeone(a_{m}(i),b_{m}(j);L(\delta,b_{m}(j))\typethree)$.

\item If $m' < m$,
then $\xi_{m'}$ is less than both
$\Delta\typefour(e_{a_m(i)},e_{a_{m+1}(i)}\typefour)$
and $\Delta\typefour(e_{b_m(j)},e_{b_{m+1}(j)}\typefour)$.

\item \label{last_entry} 
$e_{a_m(i)}\typeone(\max L(\delta,b_m(j))\typethree) >
e_{b_m(j)}\typeone(\max L(\delta,b_m(j))\typethree)$.

\item 
$\mu\typeone(\delta,b_0(j)\typeone)$ is the restriction of $\mu\typeone(\delta,b_m(j)\typeone)$
to $L\typeone(\delta,b_0(j)\typeone)$.

\item If $m' < m$, then
$\mu\typeone(\delta,b_m(j)\typeone)(\xi_{m'}) = w$.

\end{enumerate}
Now let $n$ be given.
Pick a $\gamma_0 < \delta$ which is an upper bound for each
$L\typeone(\delta,b_n(j)\typeone)$ for $j < l$ and all
$\xi < \delta$ such that for some $m,m' \leq n$ and $j < l$
$$e_{b_m(j)}(\xi) \ne e_{b_{m'}(j)}(\xi)$$
(the later set is finite by Fact \ref{rho1_fact}).
Using elementarity of $M$ and Fact \ref{L_lim},
select an $a$ in $\Acal$ such that $a < \delta$ and
for all $i < k$ and $j < l$
$$L\typeone(\delta,b_n(j)\typeone) < \Delta\typefour(e_{a(i)},e_{a_{n}(i)}\typefour),$$
$$\gamma_0 < L(a(i),\delta).$$
Now let $i < k$, $j < l$, and $m \leq n$ be fixed.
It follows from Fact \ref{L_trans} that
$$\mu\typeone(a(i),b_m(j)\typeone) = \mu\typeone(a(i),\delta\typeone)
\cup \mu\typeone(\delta,b_m(j)\typeone).$$
Also,
$$\mu\typetwo(a(i),b_m(j)\typetwo)(\xi_{m'}) =
\mu\typetwo(\delta,b_m(j)\typetwo)(\xi_{m'}) = w.$$
Finally, $e_{b_m(j)} \restriction L(a(i),\delta)$ does not
depend on $m$ and therefore
$$\Osc\typeone(a(i),b_0(j);L(a(i),\delta)\typeone) =
  \Osc\typeone(a(i),b_m(j);L(a(i),\delta)\typeone).$$
Hence, by item \ref{last_entry},
$$\Osc\typeone(a(i),b_m(j)\typeone) = \Osc\typeone(a(i),b_0(j)\typeone) \cup
\{\xi_{m'} : m' < m\}.$$
This completes the proof of Lemma \ref{osc_lem}.
\end{proof}

\section{A negative partition relation and the coloring $o$}
\label{o_section}

We will begin this section with an important consequence of
Lemma \ref{osc_lem}.
First I will introduce some notation.
\begin{defn}
The function $*:\omega \to \omega$ is defined by letting
$*(0) = 0$ and if $m > 0$, then
$*(m) = n$ if the $n\Th$ prime is the least prime
which does not divide $m$.
If $f$ is a function taking values in $\omega$, 
then $f^*$ will denote the composition of $f$ followed by $*$. 
\end{defn}
The notation for $*$ is chosen to mimic its usage in \cite{cseq}.
In particular, the following is easily verified.
\begin{fact}
If $X \subseteq \omega$ contains arbitrarily long intervals
of integers, then the image of $X$ under $*$ is $\omega$.
\end{fact}

Hence Lemma \ref{osc_lem} immediately yields the following result.

\begin{thm}
If $A,B \subseteq \omega_1$ are uncountable and $n < \omega$, then
there are $\alpha$ in $A$ and $\beta$ in $B$ with $\alpha < \beta$ and
$\osc^*(\alpha,\beta) = n$.
\end{thm}

In other words, $\osc^*$ witnesses
$\omega_1 \not \rightarrow [\omega_1;\omega_1]^2_\omega$.
If $\iota_\beta:\beta \to \omega$ is an injection for
each $\beta < \omega_1$, then for $\alpha < \beta$ we can put
$f (\alpha,\beta) = \xi$ if
$\osc^*(\alpha,\beta) = \iota_\beta(\xi)$
and $f (\alpha,\beta) = 0$ if no such $\xi$ exists.
It is easily checked this witnesses the
negative partition relation in the following theorem.

\begin{thm}
$\omega_1 \not \rightarrow [\omega_1;\omega_1]^2_{\omega_1}$.
\end{thm}

Two variations of $\omega_1 \rightarrow (\omega_1;\omega_1)^2_2$
which are closely related to (S) and (L) respectively are:
$$\omega_1 \rightarrow \typeone(\omega_1,(\omega_1;\fin)\typeone)^2,$$
$$\omega_1 \rightarrow \typeone(\omega_1,(\fin;\omega_1)\typeone)^2.$$
The statement
$\omega_1 \rightarrow \typeone(\omega_1,(\omega_1;\fin)\typeone)^2$ asserts
that if $c:[\omega_1]^2 \to 2$ then either
\begin{enumerate}

\item there is an uncountable $X \subseteq \omega_1$ such that
$c$ is constantly $0$ on $[X]^2$ or

\item there are uncountable $A \subseteq \omega_1$ and
uncountable pairwise disjoint
$\Bcal \subseteq [\omega_1]^{< \aleph_0}$ such that for all
$\alpha$ in $A$ and $b$ in $\Bcal$ with $\alpha < b$, there is a
$\beta$ in $b$ with $c(\alpha,\beta) = 1$.

\end{enumerate}
The statement
$\omega_1 \rightarrow \typeone(\omega_1,(\fin;\omega_1)\typeone)^2$ is similarly
defined by replacing $A \subseteq \omega_1$ by a pairwise disjoint
$\Acal \subseteq [\omega_1]^{<\aleph_0}$ and $\Bcal$ by
$B \subseteq \omega_1$.
These imply (S) and (L) respectively (see \cite{forcing_partition})
and are each consequences of
$\omega_1 \rightarrow (\omega_1;\omega_1)^2_2$.\footnote{
Of course this latter statement is trivial in light of Theorem
\ref{no_rectangle} but this implication was already observed by the
1980's.}

In \cite{forcing_partition}, Todor\v{c}evi\'{c} showed that
$\omega_1 \rightarrow \typeone(\omega_1,(\omega_1;\fin)\typeone)^2$ was relatively consistent
with $\ZFC$ and in fact follows from $\PFA$. 
We will now consider a map $o:[\omega_1]^2 \to \omega$
such that $o^*$ witnesses
a strong failure of
$\omega_1 \rightarrow \typeone(\omega_1,(\fin;\omega_1)\typeone)^2$.

\begin{defn}
Let $o(\alpha,\beta)$ denote
$$\sum_{q \ne 0} (|\mu(\alpha,\beta;\alpha)^{-1}(q)|\ \mod q).$$
\end{defn}

\begin{thm}\label{o_thm}
Let $\Acal \subseteq [\omega_1]^k$ and
$\Bcal \subseteq [\omega_1]^l$ be uncountable and pairwise
disjoint families.
For every $\chi:k \to 2$ and
$\pi:k \to l$ there are $a$ in $\Acal$
and $b$ in $\Bcal$ such that $a < b$ and for all $i < k$,
$$o^*\typethree(a(i),b(\pi(i))\typethree) = \chi(i).$$
\end{thm}

\begin{remark}
It is possible to prove this result if $2$ is replaced by
any positive integer.
We will only need this formulation and, alas, we are running out of
good letters for integer variables.
\end{remark}
\begin{proof}
Let $\Acal$, $\Bcal$, $\chi$, and $\pi$
be given as in the statement of the theorem.
Select distinct primes
$q_i$ $(i < k)$ which are each larger than 8.
By refining $\Acal$ if necessary, we may assume that
there is a $w$ in $C(2^{\omega},\omega)$ such that if
$a$ is in $\Acal$, then $w(a(i)) = q_i$ for all $i < k$.
Applying Lemma \ref{osc_lem}, find
$a \in \Acal$, $b_m \in \Bcal$, $\xi_m$ for each
$m$ less than $n = 6 \cdot {\displaystyle \prod_{i < k}} q_i$ such that
the following conditions are satisfied
for all $m < n$, $i < k$, and
$j < l$:
\begin{enumerate}

\item $a < b_m$.

\item $\Osc\typeone(a(i),b_m(j)\typeone)$ is the disjoint union of the sets
$\Osc\typeone(a(i),b_0(j)\typeone)$ and $\{\xi_{m'}:m'< m\}$.

\item
$\mu\typeone(a(i),b_0(j)\typeone)$ is
$\mu\typeone(a(i),b_m(j)\typeone)$
restricted to
$L\typeone(a(i),b_0(j)\typeone)$.

\item
$\mu\typeone(a(i),b_m(j)\typeone)(\xi_{m'}) = w$ whenever $m' < m$.

\end{enumerate}

Put $$r_i = \Big(\sum_{s \ne 0,q_i}
|\mu\typethree(a(i),b_0(\pi(i));a(i)\typethree)^{-1}(s)| \ 
\mod s \Big) \mod 6.$$
Use the Chinese Remainder Theorem to find an $x$
which satisfies the following equations for each $i < k$:
$$x + |\mu\typethree(a(i),b_0(\pi(i));a(i)\typethree)^{-1}(q_i)| \equiv 6 - r_i
+ 2^{\chi(i)} \qquad (\mod q_i).$$
Notice that the right hand side of these equations are always
at least 0 and strictly less than $q_i$ since $8 < q_i$.
Let $b = b_x$.

We now need to show that for each $i < k$,
the least prime which does
not divide $o\typethree(a(i),b(\pi(i))\typethree)$ is $2$ if $\chi(i) = 0$ and
$3$ otherwise.
Let $i < k$ be fixed.
Notice that, for some $y$
\begin{displaymath}
\begin{array}{rl}
o\typethree(a(i),b(\pi(i))\typethree) = & \Big(x +
\mu\typethree(a(i),b_0(\pi(i));a(i)\typethree)^{-1}(q_i)\Big) \mod q_i + r_i + 6 y \\
= & (6 - r_i + 2^{\chi(i)})
+ r_i + 6 y \\
= & 2^{\chi(i)} + 6 (y + 1). \\
\end{array}
\end{displaymath}
Since $3$ divides $6(y+1)$ but not $2^{\chi(i)}$, $3$ does not
divide $o\typethree(a(i),b(\pi(i))\typethree)$.
On the other hand, $2$ divides $o\typethree(a(i),b(\pi(i))\typethree)$ if and only
of $\chi(i) = 1$.
It follows that $o^*\typethree(a(i),b(\pi(i))\typethree) = \chi(i)$.
\end{proof}

Notice that we can think of the above result in the following way.
We say that $\chi$ is a \emph{pattern} if
for some positive integers $k$ and $l$ and some subset $D \subseteq k
\times l$, $\chi$ is a function from $D$ into $\omega$.
A coloring $c:[\omega_1]^2 \to \omega$ is said to always realize
a pattern $\chi$ if whenever
$\Acal \subseteq [\omega_1]^{k}$ and
$\Bcal \subseteq [\omega_1]^{l}$ are pairwise disjoint and
uncountable,
there are $a$ in $\Acal$ and $b$ in $\Bcal$ such that $a < b$ and
$$c\typetwo(a(i),b(j)\typetwo) = \chi(i,j)$$
whenever $(i,j)$ is in $D$.

Theorem \ref{o_thm} just says that the coloring $o^*$ always realizes
every binary pattern in which $D$ defines a function.
On the other hand, Todor\v{c}evi\'{c} has shown that
PFA implies
$\omega_1 \rightarrow \typeone(\omega_1,(\omega_1;\fin)\typeone)^2$
which implies that for every coloring $c:[\omega_1]^2 \to 2$,
there is a constant pattern on $1 \times l$
which is not always realized by $c$. 
This suggests the following question.
\begin{question}
($\PFA$)
If $c:[\omega_1]^2 \to 2$, is there a pattern on $1 \times 2$
which is not always realized?
\end{question}
The proof of Theorem \ref{X^2_discrete} below shows that $o$ cannot provide
a counterexample to this question.

\section{A family of binary relations on $\omega_1$}
\label{Tukey_section}
In this section we will consider the family of binary relations on
$\omega_1$ with the Tukey order.
Recall that the definition of the Tukey order $\leq$ on
binary relations $R$ and $S$.
\begin{defn}
If $R$ and $S$ are binary relations,
then we write $R \leq S$ iff there are functions $f:\dom(R) \to \dom(S)$ and
$g:\ran(S) \to \ran(R)$ such that
$$f(x)\ S\ y \textrm{ implies } x\ R\ g(y).$$
\end{defn}
This order was first considered by J.W. Tukey in the class of
transitive relations \cite{conv_unif_topology}.
It make sense, however, to consider this relation on the more general
class of binary relations (see \cite{Vojtas}).

Our focus will be on binary relations on $\omega_1$.
Two fundamental examples are the well order
$\omega_1$ and the family $[\omega_1]^{<\aleph_0}$ of finite
subsets of $\omega_1$ ordered by inclusion.
We will also need a few standard operations on relations.
If $R$ and $S$ are relations, then $R \oplus S$ is the relation
which is the disjoint union of $R$ and $S$.
If $m$ is a cardinal and $R$ is a relation, then we will let $m\cdot
R$ will denote the direct sum of $m$ copies of $R$.
If $R$ and $S$ are two binary relations, then
$R \meet S$ and $R \join S$ are the relations with domains
$\dom(R) \times \dom(S)$ and ranges $\ran(R) \times \ran(S)$
such that
$$(a,b)\ R \meet S\ (c,d) \textrm{ iff }
a\ R\ c \textrm{ or } b\ S\ d,$$
$$(a,b)\ R \join S\ (c,d) \textrm{ iff }
a\ R\ c \textrm{ and } b\ S\ d.$$
The join $R \join S$ is often denoted $R \times S$.
It is easily verified that $\meet$ and $\join$ give
lower and upper bounds respectively.\footnote{In general these
need not be optimal --- they are in the class of directed relations.}

It was observed that the binary relations on $\omega_1$
which one knew how to construct in
$\ZFC$ were either below $\aleph_0 \cdot \omega_1$ or above
$[\omega_1]^{<\aleph_0}$ in the Tukey order.
\begin{example}
Define $\alpha R \beta$ iff
$\alpha < \beta$ and $z_\alpha < z_\beta$.
\end{example}
The relation $R$ is essentially Sierpi\'nski's partition which witnesses
$\omega_1 \not \rightarrow (\omega_1)^2_2$.
In this case $R \leq \aleph_0 \cdot \omega_1$ by the following
Tukey maps:
$$f(\alpha) = (\min\{n:z_\alpha < q_n\},\alpha)$$
$$g(\alpha,n) = \min\{\beta:\alpha < \beta \mand q_n < z_\beta\}$$
where $\{q_n\}_{n < \omega}$ is an enumeration of a countable order-dense
subset of $2^\omega$.
On the other hand the following example shows that with an additional
hypothesis such as $\diamond$, one can construct more complex relations.
\begin{example}
Let $R$ be a tree order on a set $T$.
If $R \leq \aleph_0 \cdot \omega_1$, then
$T$ is the union of countably many $R$-chains.
If $[\omega_1]^{<\aleph_0} \leq R$, then
$T$ contains an uncountable $R$-antichain.
In particular, if $(T,R)$ is a Suslin tree, then
$R \not \leq \aleph_0 \cdot \omega_1$ and
$[\omega_1]^{< \aleph_0} \not \leq R$.
\end{example}
These observations led to the following conjecture.
\begin{conjecture-section}\label{false_conjecture} ($\PFA$)
If $R$ is a binary relation, then either
$R \leq \aleph_0 \cdot \omega_1$ or else $[\omega_1]^{< \aleph_0} \leq R$.
\end{conjecture-section}
This can be considered a basis conjecture in the following sense.
Let $\Rcal$ denote the class of all binary relations $R$ on $\omega_1$
such that $R$ is not reducible to $\aleph_0 \cdot \omega_1$.
Then the above conjecture is just the assertion that $\Rcal$ has
a single element basis consisting of $[\omega_1]^{<\aleph_0}$.
This conjecture was given further plausibility by the following
theorem which implies that it is true for transitive relations.
\begin{thm} \cite{class_trans}
($\PFA$)
Every transitive relation on $\omega_1$ is Tukey equivalent to
one of the following for some non negative integers $n_i$ $(i < 5)$:
\begin{enumerate}

\item $n_0 \cdot 1 \oplus n_1 \cdot \omega \oplus n_2 \cdot \omega_1
\oplus n_3 \cdot \omega \times \omega_1 \oplus n_4
\cdot [\omega_1]^{<\aleph_0}$.

\item $\aleph_0 \cdot 1 \oplus n_2 \cdot \omega_1 \oplus n_3 \cdot
\omega \times \omega_1 \oplus n_4 \cdot \fin$.

\item $\aleph_0 \cdot \omega_1 \oplus n_4 \cdot \fin$.

\item $\aleph_0 \cdot \fin$.

\item $=$.

\end{enumerate} 
\end{thm}
We will now see, however, that Conjecture \ref{false_conjecture}
is provably false.
\begin{defn}
Let $c(\alpha,\beta)$ denote $o(\alpha,\beta) \mod 2$ if
$\alpha < \beta$.
\end{defn}
\begin{defn}
Let $R$ denote the relation with domain and range $\omega_1$
defined by letting $\alpha R \beta$ iff
$\alpha = \beta$ or
$\alpha < \beta$ and $c(\alpha,\beta) = 1$.
If $X$ is an uncountable subset of $\omega_1$, let $R_X$
denote the relation with
domain $X$ and range $\omega_1$ which is the
restriction of $R$ to $X \times \omega_1$.
\end{defn}
Trivially, $\omega_1 \leq R_X$ for every uncountable
$X$ and if $X \subseteq Y$,
then $R_X \leq R_Y$.
The following shows that none of these relations are above
$[\omega_1]^{< \aleph_0}$ in the Tukey order.

\begin{thm}
$[\omega_1]^{<\aleph_0} \not \leq R$.
\end{thm}

\begin{proof}
Suppose that $f:[\omega_1]^{<\aleph_0} \to \omega_1$ and
$g:\omega_1 \to [\omega_1]^{< \aleph_0}$ are functions.
It suffices to show that they are not Tukey reductions.

Applying the pressing down lemma,
find a stationary $B \subseteq \omega_1$ and a finite
$x_0 \subseteq \omega_1$ such that for all $\beta$ in $B$
$$g(\beta) \cap \beta = x_0.$$
Let
$$\Xcal = \{x \in [\omega_1]^{< \aleph_0}:x \not \subseteq x_0\},$$
$$A = \{f(x):x \in \Xcal\}.$$
By refining $B$ if necessary, we may assume that if $\alpha$ is in $A$
and $\beta$ is in $B$ with $\alpha < \beta$, then there is an $x$ in 
$\Xcal$ such that $f(x) =\alpha$ and $x \subseteq \beta$.

If $A$ is countable, then we can find a $\alpha$ in $A$ such that
$$\{x \in \Xcal:f(x) = \alpha\}$$ is uncountable.
It is then easy for find an $x$ in $\Xcal$ such that
$x$ is not contained in $g(\alpha)$ but $f(x) = \alpha$ and hence
$f(x)\ R\ \beta$, witnessing that $f,g$ are not Tukey maps.

If $A$ is uncountable, apply Theorem \ref{o_thm} to obtain
an $\alpha$ in $A$ and a $\beta$ in $B$ such that
$\alpha < \beta$ and $c(\alpha,\beta) = 1$.
Now pick an $x$ in $\Xcal$ such that $f(x) = \alpha$ and $x \subseteq \beta$
and observe that $x$ is not contained in $g(\beta)$ and
yet $f(x)\ R\ \beta$.
\end{proof}

\begin{thm}
For all finite sets $F$ of uncountable subsets of $\omega_1$,
$$\bigmeet_{X \in F} R_X \not \leq \aleph_0 \cdot \omega_1.$$
\end{thm}

\begin{proof}
By replacing $F$ with a disjoint refinement if necessary, we may
assume without loss of generality that $F$ consists of pairwise
disjoint sets.
Suppose for contradiction that
$$f:\prod_{X \in F} X \to \omega \times \omega_1$$
$$g:\omega \times \omega_1 \to \prod_{X \in F} \omega_1$$
are Tukey maps.
Find  uncountable pairwise disjoint
$\Acal, \Bcal \subseteq [\omega_1]^{|F|}$ and an $n$ such that:
\begin{enumerate}

\item For each element $a$ of $\Acal \cup \Bcal$,
$a \cap X$ is a singleton for all $X$ in $F$.

\item
$f \restriction \Acal$ is an injection into $\{n\}\times \omega_1$.

\item
$g^{-1} \Bcal \subseteq \{n\} \times \omega_1$.

\end{enumerate}
Applying Theorem \ref{o_thm}, it is possible to find an $a$ in $\Acal$
and an $x$ in $g^{-1} \Bcal$ such that $a < g(x)$,
$f(a) < x$, and for all $X$ in $F$
$$c(a \cap X,g(x) \cap X) = 0$$
(where we are identifying the singletons
$a \cap X$ and $g(x) \cap X$ with their element).
It follows immediately that
$f(a) \leq x$
while $a$ is not $\bigmeet_{X \in F} R_X$-related to $g(x)$, a contradiction. 
\end{proof}

The next theorem show that these relations are typically
incomparable. 

\begin{thm}
If $X \setminus Y$ is uncountable, then
$R_X \not \leq R_Y$.
\end{thm}

\begin{proof}
Suppose for contradiction that $f:X \to Y$ and $g:\omega_1 \to \omega_1$
are Tukey maps.
To obtains a contradiction,
it suffices to find $\alpha$ in $X$ and a $\beta$ in $\omega_1$
such that
$f(\alpha) < \beta$,
$\alpha < g(\beta)$,
$c(f(\alpha),\beta) =1$, and
$c\typeone(\alpha,g(\beta)\typeone) =0$.
Since $f$ and $g$ are Tukey reductions it is possible to
find a uncountable $X_0 \subseteq X \setminus Y$ and
uncountable $Z \subseteq \omega_1$ such that $f \restriction X_0$
and $g \restriction Z_0$ are injections into sets disjoint from
$X_0$ and $Z_0$, respectively, and the inequalities
$\alpha < f(\alpha)$, $\beta = f(\beta)$, $\beta < f(\beta)$
are uniformly true or false as $\alpha$ ranges over
$X_0$ and $\beta$ ranges over $Z_0$.
Put
$$\Acal= \{\{\alpha,f(\alpha)\}:\alpha \in X_0\}$$
$$\Bcal= \{\{\beta,g(\beta)\}:\beta \in Z_0\}$$
and notice that these families are uncountable and consist of
pairwise disjoint sets of uniform cardinality which, in the case of
$\Acal$, is two.
By our uniformity assumption on the relations
$<$ and $=$, it is possible to find $i_0$ and $j_0$ such that
if $a$ is in $\Acal$, then $a(i_0)$ is in $X_0$ and
if $b$ is in $\Bcal$, then $b(j_0)$ is in $Z_0$.
Put $\pi(i_0) = j_0$ and $\pi(1-i_0) = 1-j_0$ if
$\Bcal$ consists of doubletons and $\pi(1-i_0) = j_0$
otherwise.
Define $\chi(i_0) = 0$ and $\chi(1-i_0) = 1$.
Applying Theorem \ref{o_thm} to $\Acal$, $\Bcal$, $\pi$, and $\chi$,
there are $a$ in $\Acal$ and $b$ in $\Bcal$ such that
$$a < b$$
$$c\typeone(a(i),b(\pi(i))\typethree) = \chi(i).$$
Translating the outcome, it is easily
checked that $\alpha = a(i_0)$ and $\beta = b(j_0)$ are as desired.
\end{proof}

The following remain unclear.
\begin{question}
For which families $\Fcal \subseteq [\omega_1]^{\aleph_1}$ is there
an $S$ such that $S \not \leq \aleph_0 \cdot \omega_1$ and
$S \leq R_X$ for all $X$ in $\Fcal$.
\end{question}

\begin{question}
Is the collection of all $S$ with
$S \not \leq \aleph_0 \cdot \omega_1$
downwards-directed in the Tukey order?
\end{question}

\section{An L space and the non-existence of a small
basis for the regular Hausdorff spaces.}
\label{Lspace_section}
In this section I will give an example of an \emph{L space} ---
a hereditarily Lindel\"of, not hereditarily
separable space.\footnote{According
to Juhasz, M.E. Rudin coined the phrase ``S space'' but to mean
an hereditarily separable space.
Juhasz then introduced the phrase ``L space''
and gave the current definitions to S and L spaces as those
which distinguish hereditary separability from the hereditary
Lindel\"of property.
Also, it was shown in \cite{alphaS_alphaL} that examples exist
without the assumption that the spaces are regular and Hausdorff.}
The question of the existence of such spaces was first asked explicitly
in \cite{alphaS_alphaL}, though arguably this question can be traced
to Suslin's \cite{Suslin_problem} where he posed his famous hypothesis.
For instance, an immediate consequence of
Kurepa's \cite{orders_ramifications} is that a Suslin line
is an example of an L space.
From the 1960's until the 1980's, there was a concerted effort
to understand both L spaces and their dual the S
space.
I refer the reader to Juhasz's \cite{SL_spaces:Juhasz} and Rudin's
\cite{SL_spaces:Rudin} as well as
Roitman's more recent \cite{basic_SL} for more discussion on these
developments.
I will mention a few.
First, Roitman and Zenor showed that there was a relationship between
the existence of certain S and L spaces.
\begin{thm} \cite{cohen_random} \cite{HmS_HmL}
There is a strong S space iff there is a strong L space.\footnote{Here
a space is a strong S (L) space iff all of its finite powers are S (L) spaces.}
\end{thm}
Hence the difference in the existence of S and L spaces lies in the
properties of their finite powers.
This gives some explanation as to why the existence of S and L spaces
seem to be such similar hypotheses at first.

There are a number of results under $\MA_{\aleph_1}$ which limit the
existence of S and L spaces.
\begin{thm} \cite{strong_S_L} ($\MA_{\aleph_1}$)
There are no strong S or L spaces.
\end{thm}

\begin{thm} \cite{S_L_MA} ($\MA_{\aleph_1}$)\label{no_1stctblL}
There are no first countable L spaces.
\end{thm}

\begin{thm} \cite{consequences_MA} ($\MA_{\aleph_1}$)
If $K$ is a compact space which contains an L space, then
$K$ maps continuously onto $[0,1]^{\aleph_1}$. 
\end{thm}

Finally, Todor\v{c}evi\'{c} proved that S spaces do not exist
assuming $\PFA$.
\begin{thm} \cite{forcing_partition} ($\PFA$)
If a space is hereditarily separable, then it is hereditarily
Lindel\"of.
\end{thm}

Until now, however, it remained unclear whether Todor\v{c}evi\'c's methods
could be used to prove that $\PFA$ implies there are no L spaces.
I will now show that this is not the case --- that there is an L space
which can be constructed without appealing to additional axioms of set theory.
For each $\alpha < \omega_1$, put
$$W_\alpha = \{\alpha\} \cup \{\beta > \alpha: c(\alpha,\beta) =
1\}.$$
If $X \subseteq \omega_1$, then let
$\tau[X]$ be the topology on $X$ generated by declaring
$W_\xi \cap X$ to be clopen for all $\xi$ in $X$.
If $X$ is uncountable, then clearly $\tau[X]$ is non-separable.
The following theorem shows that $\tau[X]$ is always Fr\'echet.

\begin{thm}
If $X$ is a set and $\Fcal$ is a point
countable point separating family of subsets of $X$,
then the topology on $X$ defined by
declaring elements of $\Fcal$ to be clopen is countably tight and
has every countable subspace metrizable.
In particular the topology is Fr\'echet.
\end{thm}

\begin{remark}
I attribute this result to Zoltan Balogh.
He once told me that ``If there is an L space, then there is a countably
tight one.''
While I never saw his proof, I think it is reasonable to assume this
he may have proceeded along these lines.
\end{remark}

\begin{proof}
That countable subspaces are metrizable follows from the assumptions and
the well known fact that a (regular Hausdorff)
second countable spaces are metrizable.
To see that the space is countably tight, let
$A \subseteq X$ have an accumulation point $x$ in $X$.
Let $M$ be a countable elementary submodel of $H(\theta)$ for $\theta$ regular
and large enough so that $X$, $\Fcal$, $x$, and $A$ are all in $M$.
It suffices to show that $x$ is an accumulation point of $A \cap M$.

Suppose this is not the case and let
$U_i$ $(i < k)$ and $V_i$ $(i < l)$ be elements of $\Fcal$ such that
$$W = \bigcap_{i < k} U_i \setminus \bigcup_{j < l} V_j$$
contains $x$ and is disjoint from $A \cap M$.
The important observation is this:
if $V$ is in $\Fcal$, then $V \cap M$ is non-empty
iff $V$ is in $M$.
This is because $\{y \in V:V \in \Fcal\}$ is countable for all
$y$ in $X$ and therefore a subset of $M$,
\ref{definable_is_in}) whenever
$y$ is in $M$ (by Facts \ref{uncountable_test1}.
Hence each $U_i$ must be in $M$ since $x$ is in $U_i$ and $x$ is in $M$.
Moreover, since we are only interested in having $W$ be disjoint from
$A \cap M$, we may assume without loss of generality that each $V_j$ is
in $M$.
But then $W$ must be in $M$ and, by elementarity of $M$, there must be an
element of $W \cap A$ which is in $M$,
a contradiction.
\end{proof}

It is straightforward to use Theorem \ref{o_thm}
to show that every $\tau[X]$ is hereditarily Lindel\"of;
I will prove a result of independent interest and then derive this as
a consequence.
\begin{thm} \label{rigid}
If $X, Y \subseteq \omega_1$
have countable intersection, then
there is no continuous injection from any uncountable subspace of
$(X,\tau[X])$ into $(Y,\tau[Y])$.
\end{thm}

\begin{proof}
Without loss of generality $X$ and $Y$ are actually disjoint.
Suppose for contradiction that such an injection
$f:X_0 \to Y$ does exist
for some uncountable $X_0 \subseteq X$.
Let $\Bcal$ be the set of all $\{\beta,f(\beta)\}$ such that
$\beta$ is in $X_0$.
For simplicity we will assume that $\beta < f(\beta)$ for all $\beta$;
a similar argument is given if $f(\beta) < \beta$ for uncountably
many $\beta$.
For each $\beta$ in $X_0$, pick a basic open set $U_\beta$ which contains
$\beta$ such that $f^{-1} W_{f(\beta)}$ contains $U_{\beta}$.
By refining $X_0$ if necessary, we may assume that for some $k$
and $\chi_0:k \to 2$ there
is a $\Delta$-system $\{F_\beta:\beta \in X_0\}$ with root $F$ of
$|F| + k$-element subsets of $X$ such that:
\begin{enumerate}

\item
If $F_\alpha < \beta$, then $\beta$ is in
$U_\alpha$ iff for all $i < |F| + k$,
$c(F_\alpha(i),\beta) = \chi_0(i)$.

\item
$l = |F_\alpha \cap f(\alpha)|$ does not depend on $\alpha$.

\item $\max F$ is less than both $\alpha$
and $\min(F_\alpha \setminus F)$.

\end{enumerate} 
Let $\Acal$ be the collection of all sets of the form
$F_\alpha \setminus F \cup \{f(\alpha)\}$.
Notice that since $X$ is disjoint from $Y$, 
$f(\alpha)$ is not in $F_\alpha$ for any $\alpha$ in $X_0$.
Define $\chi:k+1 \to 2$ by $\chi(i) = \chi_0(|F| + i)$ if
$i < l$, $\chi(l) = 0$, and $\chi(i) = \chi_0(|F| + i-1)$ if $i > l$.
Define $\pi:k+1 \to 2$ by $\pi(i) = 0$ if $i \ne l$ and
$\pi(l) = 1$.
Applying Theorem \ref{o_thm}, there are $a \in \Acal$ and $b \in \Bcal$
such that $a < b$ and for all $i < k+1$
$$c\typethree(a(i),b(\pi(i))\typethree) = \chi(i).$$
Let $\alpha,\beta$ in $X_0$ be such that
$$a = (F_\alpha \setminus F) \cup \{f(\alpha)\}$$
$$b = \{\beta,f(\beta)\}.$$

To derive a contradiction, it suffices to show that $\beta$ is in
$U_\alpha$ but that $f(\beta)$ is not in $W_{f(\alpha)}$.
The latter is just a reformulation of
$$c\typethree(a(l),b(\pi(l))\typethree) = c\typeone(f(\alpha),f(\beta)\typeone) = \chi(l) = 0.$$
For the former, we need to show that for all $i < |F| + k$,
$$c(F_\alpha(i),\beta) = \chi_0(i).$$
If $i < |F|$, then $F_\alpha(i) = F_\beta(i)$
and, since $\beta$ is in $U_\beta$ and $F < \beta$,
$$c(F_\alpha(i),\beta)= c(F_\beta(i),\beta) = \chi_0(i) = 1.$$
If $i \geq |F|$, then
$\chi_0(i) = \chi(i - |F|)$ if $i - |F| < l$ and
$\chi_0(i) = \chi(i - |F|+1)$ if $i - |F| \geq l$.
If $i - |F| < l$, then
$$c(F_\alpha(i),\beta) = c(a(i-|F|),\beta) = \chi(i - |F|) = \chi_0(i)$$
and if $i - |F| \geq l$, then
$$c(F_\alpha(i),\beta) = c(a(i-|F|+1),\beta) = \chi(i - |F|+1) = \chi_0(i).$$
This finishes the proof.
\end{proof}

\begin{cor}
For every $X$, $\tau[X]$ is hereditarily Lindel\"of.
\end{cor}

\begin{proof}
If not, then $(X,\tau[X])$ would contain an uncountable discrete subspace.
It is then possible to find disjoint $Y,Z \subseteq X$ such that
$(Y,\tau[Y])$ and $(Z,\tau[Z])$ each contain uncountable discrete subspaces.
But any function from a discrete space to another discrete space is
continuous and we can then easily contradict Theorem \ref{rigid}.
\end{proof}

Kunen has shown, however, that under $\MA_{\aleph_1}$ every L space contains
an uncountable discrete subspace in one of its finite powers \cite{strong_S_L}.
As might be expected, this happens at the first possible instance in
our example.

\begin{thm}\label{X^2_discrete}
For every uncountable $X \subseteq \omega_1$,
$\tau[X]^2$ contains an uncountable discrete subspace.
\end{thm}

\begin{proof}
This essentially follows from the following observation which is of
independent interest.
\begin{prop} \label{o_Aronszajn}
The tree
$$T(o) = \{o(\cdot,\beta) \restriction \alpha:\alpha \leq \beta < \omega_1\}$$
is Aronszajn.
\end{prop}

\begin{proof}
Theorem \ref{o_thm} implies that $T(o)$ has no uncountable branches.
It suffices to prove that all levels of $T(o)$ are countable.
For this, it is sufficient to show that
$T(\rho_1)$ and $T(\mu)$ have countable levels.
In the case of $T(\rho_1)$, this is a trivial consequence
of Fact \ref{rho1_fact}.
In order to see that $T(\mu)$ has countable levels, let
$\alpha < \omega_1$ be given.
Let $\beta < \omega_1$ be greater than $\alpha$.
Using the compactness of $\alpha + 1$, select a finite set
$F_\beta \subseteq \alpha + 1$ containing $0$ and $\alpha$
so that if $\gamma_0 < \gamma$ are consecutive elements of $F_\beta$,
then
$$L(\gamma,\beta) < L(\xi,\gamma)$$
whenever $\gamma_0 < \xi < \gamma$.
It suffices to show that if
$$\mu(\cdot,\beta) \restriction F_\beta =
\mu(\cdot,\beta') \restriction F_{\beta'},$$
then
$$\mu(\cdot,\beta) \restriction \alpha = L(\cdot,\beta') \restriction
\alpha.$$
To see this, let $\xi < \alpha$ be arbitrary
and pick $\gamma_0 < \gamma$ in $F_\beta$
so that $\gamma_0 < \xi \leq \gamma$.
If $\xi = \gamma$, then $\mu(\xi,\beta) = \mu(\xi,\beta')$ by
assumption.
Otherwise,
$$\mu(\xi,\beta) = \mu(\xi,\gamma)\cup\mu(\gamma,\beta)
= \mu(\xi,\gamma)\cup\mu(\gamma,\beta') = 
\mu(\xi,\beta').$$
\end{proof}

To finish the proof of Theorem \ref{X^2_discrete},
select a sequence $(\beta_\xi^0,\beta_\xi^1)$ indexes by an
uncountable set $\Xi \subseteq \omega_1$
such that for each $\xi$ in $\Xi$ the following conditions are satisfied:
\begin{enumerate}

\item $\beta_\xi^0 < \beta_\xi^1$ are both elements of $X$.

\item Both $\beta_\xi^0$ and $\beta_\xi^1$ are greater than $\xi$.

\item If $\eta < \xi$ is in $\Xi$, then $\beta_\eta^1 < \xi$.

\item
$o(\cdot,\beta_\xi^0) \restriction \xi = o(\beta_\xi^1) \restriction \xi$.

\end{enumerate}
This is possible since $T(o)$ is Aronszajn.
Now consider the clopen neighborhoods
$$U_\xi =
(W_{\beta_\xi^1} \setminus W_{\beta_\xi^1}) \times W_{\beta_\xi^1}$$
of $(\beta_\xi^0,\beta_\xi^1)$ for $\xi$ in $\Xi$.
If $\eta < \xi$ are in $\Xi$, then
$(\beta_\eta^0,\beta_\eta^1)$ is not in $U_\xi$.
If $\xi < \eta$ are in $\Xi$, then
$$o(\beta_\xi^1,\beta_\eta^0) = o(\beta_\xi^1,\beta_\eta^1).$$
It follows that $\beta_\eta^0$ is in $W_{\beta_\xi^1}$ iff
$\beta_\beta^1$ is.
This means that $(\beta_\eta^0,\beta_\eta^1)$ is not in $U_\xi$.
Hence $(\beta_\eta^0,\beta_\eta^1)$ is in $U_\xi$ iff
$\eta = \xi$ and therefore
$$\{(\beta_\xi^0,\beta_\xi^1):\xi \in \Xi\}$$
is an uncountable discrete subspace of $\tau[X]^2$ as desired.
\end{proof}

The proof of Theorem \ref{no_1stctblL} can be used to prove
the following.
\begin{thm}
For every $X \subseteq \omega_1$,
no uncountable subspace of $\tau[X]$ is first countable.
\end{thm}

\begin{proof}
The proof of Theorem \ref{no_1stctblL} actually shows that
if an L space is not first countable, then there is a c.c.c. forcing
which destroys it.
It is easily checked that the any uncountable subspace
of $\tau[X]$ is a c.c.c. indestructible L space.
\end{proof}

We also have the following consequence of Proposition \ref{o_Aronszajn}.
\begin{thm}
For every $X$, if $f$ is a continuous function from
$\tau[X]$ into $[0,1]$, the range of $f$ is countable.
\end{thm}

\begin{proof}
For each rational $q$ in $[0,1]$, let
$U_q$ be the $f$-preimage of $[0,q)$.
Applying the hereditary Lindel\"of property of $\tau[X]$,
find countable $\Ucal_q$ consisting of basic clopen sets
such that $U_q = \cup \Ucal_q$.
Now let $\delta$ be an ordinal larger than any ordinal
mentioned in the definition of an element of some $\Ucal_q$.
Applying Proposition \ref{o_Aronszajn},
pick $\beta_0$ and $\beta_1$ in $X$
such that
$c(\cdot,\beta_0) \restriction \delta = c(\cdot,\beta_1) \restriction
\delta$.
It suffices to show that $f(\beta_0) = f(\beta_1)$.
If not, pick a rational such that $[0,q)$ contains exactly
one of $f(\beta_0)$, $f(\beta_1)$.
But then there is an element of $\Ucal_q$ which contains exactly one
of $\beta_0$, $\beta_1$.
This means that for some $\alpha < \delta$,
$$c(\alpha,\beta_0) \ne c(\alpha,\beta_1),$$
contradicting the definition of $\delta$ and our choice of
$\beta_0$, $\beta_1$.
\end{proof}

\begin{thm} \ref{no_top_basis}
Any basis for the uncountable regular Hausdorff spaces
must have cardinality at least $\aleph_2$.
\end{thm}

\begin{proof}
The first conclusion follows from Theorem \ref{rigid} and the 
observation that there is an almost disjoint
family of uncountable subsets of $\omega_1$ of cardinality $\aleph_2$.
\end{proof}

\end{document}